\documentclass[12pt]{article}
\usepackage{amsfonts}

\usepackage{amsmath}
\usepackage{epsfig}

\begin{document}

\title{\textbf{A data-reconstructed fractional volatility model}}
\author{\textbf{R. Vilela Mendes}\thanks{%
Centro de Matem\'{a}tica e Aplica\c{c}\~{o}es Fundamentais and Universidade
T\'{e}cnica de Lisboa, e-mail: vilela@cii.fc.ul.pt} and \textbf{M. J.
Oliveira}\thanks{%
Centro de Matem\'{a}tica e Aplica\c{c}\~{o}es Fundamentais and Universidade
Aberta, oliveira@cii.fc.ul.pt} \\
Complexo Interdisciplinar,\\
Av. Prof. Gama Pinto 2, 1649-003 Lisboa, Portugal}
\date{}
\maketitle

\begin{abstract}
Based on criteria of mathematical simplicity and consistency with empirical
market data, a stochastic volatility model is constructed, the volatility
process being driven by fractional noise. Price return statistics and
asymptotic behavior are derived from the model and compared with data.
Deviations from Black-Scholes and a new option pricing formula are also
obtained.
\end{abstract}


\section{Introduction}

Classical Mathematical Finance has, for a long time, been based on the
assumption that the price process of market securities may be approximated
by geometric Brownian motion 
\begin{equation}
\begin{array}{lll}
dS_{t} & = & \mu S_{t}dt+\sigma S_{t}dB\left( t\right)
\end{array}
\label{1.00}
\end{equation}
In liquid markets the autocorrelation of price changes decays to negligible
values in a few minutes, consistent with the absence of long term
statistical arbitrage. Geometric Brownian motion models this lack of memory,
although it does not reproduce the empirical leptokurtosis. On the other
hand, nonlinear functions of the returns exhibit significant positive
autocorrelation. For example, there is volatility clustering, with large
returns expected to be followed by large returns and small returns by small
returns (of either sign). This, together with the fact that autocorrelations
of volatility measures decline very slowly\cite{Ding2} \cite{Harvey} \cite
{Crato}, has the clear implication that long memory effects should somehow
be represented in the process and this is not included in the geometric
Brownian motion hypothesis.

One other hand, as pointed out by Engle\cite{Engle}, when the future is
uncertain investors are less likely to invest. Therefore uncertainty
(volatility) would have to be changing over time. The conclusion is that a
dynamical model for volatility is needed and $\sigma $ in Eq.(\ref{1.00}),
rather than being a constant, becomes a process by itself. This idea led to
many deterministic and stochastic models for the volatility (\cite{Taylor} 
\cite{Engle2} and references therein).

Using, at each step, both a criteria of mathematical simplicity and
consistency with market data, a stochastic volatility model is constructed
here, with volatility driven by fractional noise. It appears to be the
minimal model consistent both with mathematical simplicity and the market
data. It turns out that this data-inspired model is different from the many
stochastic volatility models that have been proposed in the literature. The
model will be used to compute the price return statistics and asymptotic
behavior, which are compared with actual data. Deviations from the classical
Black-Scholes result and a new option pricing formula are also obtained.

\section{The induced volatility process}

The basic hypothesis for the model construction are:

(H1) The log-price process $\log S_{t}$ belongs to a probability product
space $\Omega \otimes \Omega ^{^{\prime }}$ of which the first one, $\Omega $%
, is the Wiener space and the second, $\Omega ^{^{\prime }}$, is a
probability space to be characterized later on. Denote by $\omega \in \Omega 
$ and $\omega ^{^{\prime }}\in \Omega ^{^{\prime }}$ the elements (sample
paths) in $\Omega $ and $\Omega ^{^{\prime }}$ and by $\mathcal{F}_{t}$ and $%
\mathcal{F}_{t}^{^{\prime }}$ the $\sigma -$algebras in $\Omega $ and $%
\Omega ^{^{\prime }}$ generated by the processes up to $t$. Then, a
particular realization of the log-price process is denoted 
\[
\log S_{t}\left( \omega ,\omega ^{^{\prime }}\right) 
\]
This first hypothesis is really not limitative. Even if none of the
non-trivial stochastic features of the log-price were to be captured by
Brownian motion, that would simply mean that $S_{t}$ is a trivial function
in $\Omega $.

(H2) The second hypothesis is stronger, although natural. We will assume
that for each fixed $\omega ^{^{\prime }}$, $\log S_{t}\left( \bullet
,\omega ^{^{\prime }}\right) $ is a square integrable random variable in $%
\Omega $.

\begin{center}
---------
\end{center}

From the second hypothesis it follows that, for each fixed $\omega
^{^{\prime }}$, 
\begin{equation}
\begin{array}{lll}
\frac{dS_{t}}{S_{t}}\left( \bullet ,\omega ^{^{\prime }}\right) & = & \mu
_{t}\left( \bullet ,\omega ^{^{\prime }}\right) dt+\sigma _{t}\left( \bullet
,\omega ^{^{\prime }}\right) dB\left( t\right)
\end{array}
\label{2.1}
\end{equation}
where $\mu _{t}\left( \bullet ,\omega ^{^{\prime }}\right) $ and $\sigma
_{t}\left( \bullet ,\omega ^{^{\prime }}\right) $ are well-defined processes
in $\Omega $. (Theorem 1.1.3 in Ref.\cite{Nualart})

Recall that if $\left\{ X_{t},\mathcal{F}_{t}\right\} $ is a process such
that 
\begin{equation}
\begin{array}{lll}
dX_{t} & = & \mu _{t}dt+\sigma _{t}dB\left( t\right)
\end{array}
\label{2.2}
\end{equation}
with $\mu _{t}$ and $\sigma _{t}$ being $\mathcal{F}_{t}-$adapted processes,
then 
\begin{equation}
\begin{array}{lll}
\mu _{t} & = & \underset{\varepsilon \rightarrow 0}{\lim }\frac{1}{%
\varepsilon }\left\{ \left. E\left( X_{t+\varepsilon }-X_{t}\right) \right| 
\mathcal{F}_{t}\right\} \\ 
\sigma _{t}^{2} & = & \underset{\varepsilon \rightarrow 0}{\lim }\frac{1}{%
\varepsilon }\left\{ \left. E\left( X_{t+\varepsilon }-X_{t}\right)
^{2}\right| \mathcal{F}_{t}\right\}
\end{array}
\label{2.3}
\end{equation}

The process associated to the probability space $\Omega ^{^{\prime }}$ is
now to be inferred from the data. According to (\ref{2.3}), for each fixed $%
\omega ^{^{\prime }}$ realization in $\Omega ^{^{\prime }}$ one has 
\begin{equation}
\sigma _{t}^{2}\left( \bullet ,\omega ^{^{\prime }}\right) =\underset{%
\varepsilon \rightarrow 0}{\lim }\frac{1}{\varepsilon }\left\{ E\left( \log
S_{t+\varepsilon }-\log S_{t}\right) ^{2}\right\}  \label{2.4}
\end{equation}
Each set of market data corresponds to a particular realization $\omega
^{^{\prime }}$. Therefore, assuming the realization to be typical, the $%
\sigma _{t}^{2}$ process may be reconstructed from the data by the use of (%
\ref{2.4}). To this data-reconstructed $\sigma _{t}$ process we call the 
\textit{induced volatility}.

For practical purposes we cannot strictly use Eq.(\ref{2.4}) to reconstruct
the induced volatility process, because when the time interval $\varepsilon $
is very small the empirical evaluation of the variance becomes unreliable.
Instead, we estimate $\sigma _{t}$ from 
\begin{equation}
\sigma _{t}^{2}=\frac{1}{\left| T_{0}-T_{1}\right| }\mathnormal{var}\left(
\log S_{t}\right)  \label{2.5}
\end{equation}
with a time window $\left| T_{0}-T_{1}\right| $ sufficiently small to give a
reasonably local characterization of the volatility, but also sufficiently
large to allow for a reliable estimate of the local variance of $\log S_{t}$.

As an example, daily data has been used with time windows of 5 to 9 days.
The upper left panel of Fig.1 shows the result of application of (\ref{2.5})
to the New York Stock Exchange (NYSE) aggregate index in the period 1966$-$%
2000, with a time window $\left| T_{0}-T_{1}\right| =5$ days. Notice that to
discount trend effects and approach asymptotic stationarity of the process,
before application of (\ref{2.5}), the data has been detrended and rescaled
as explained in Ref.\cite{Vilela1}. Namely, a polynomial fit is performed
for increasing orders until the fitted polynomial is no longer well
conditioned. This seems to be a reasonable detrending method insofar as it
leads to an asymptotically stationary signal \cite{Vilela1}.

Then, as a first step towards finding a mathematical characterization of the 
\textit{induced volatility process} one looks for scaling properties. Namely
one checks whether a relation of the form 
\begin{equation}
E\left| \sigma \left( t+\Delta \right) -\sigma \left( t\right) \right| \sim
\Delta ^{H}  \label{2.6}
\end{equation}
or 
\begin{equation}
E\left| \frac{\sigma \left( t+\Delta \right) -\sigma \left( t\right) }{%
\sigma \left( t\right) }\right| \sim \Delta ^{H}  \label{2.7}
\end{equation}
holds for the induced volatility process. This would be the behavior implied
by most stochastic volatility models proposed in the past. It turns out that
the data shows this to be a very bad hypothesis, meaning that the induced
volatility process itself is not self-similar.

Instead, using a standard technique to detect long-range dependencies\cite
{Taqqu}, one computes the empirical integrated log-volatility and finds that
it is well represented by a relation of the form 
\begin{equation}
\sum_{n=0}^{t/\delta }\log \sigma \left( n\delta \right) =\beta t+R_{\sigma
}\left( t\right)  \label{2.8}
\end{equation}
where, as shown in the lower right panel of Fig.1, the $R_{\sigma }\left(
t\right) $ process has very accurate self-similar properties ($\delta =1$
day for daily data).

\begin{figure}[tbh]
\begin{center}
\psfig{figure=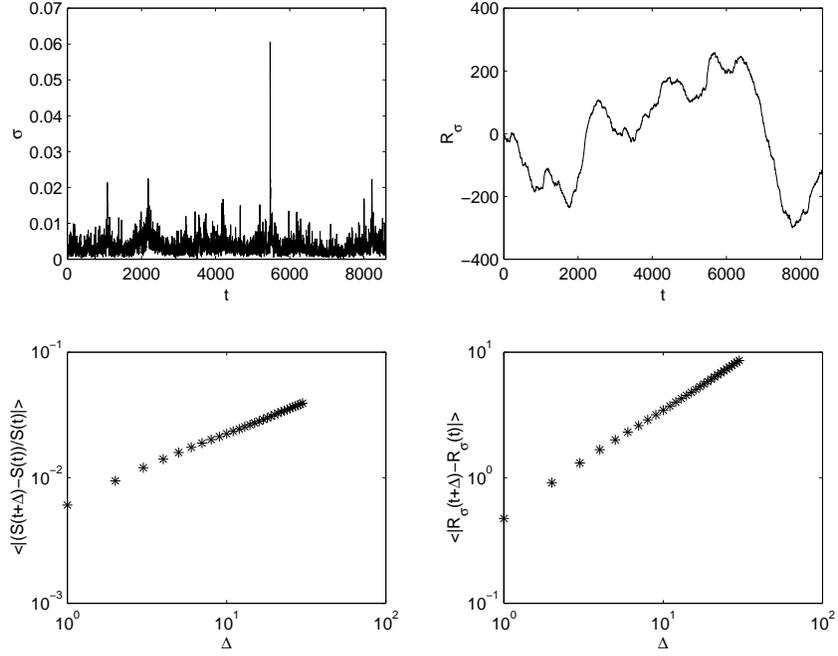,width=11truecm}
\end{center}
\caption{Self-similar properties of the integrated log-volatility $\beta
t+R_{\sigma }\left( t\right) $ process}
\end{figure}

This suggests the following mathematical identification:

(a) Recall that if a nondegenerate process $X_{t}$ has finite variance,
stationary increments and is self-similar 
\begin{equation}
\mathnormal{Law}\left( X_{at}\right) =\mathnormal{Law}\left(
a^{H}X_{t}\right)  \label{2.9}
\end{equation}
then \cite{Embrechts} $0<H\leq 1$ and 
\begin{equation}
\mathnormal{Cov}\left( X_{s},X_{t}\right) =\frac{1}{2}\left( \left| s\right|
^{2H}+\left| t\right| ^{2H}-\left| s-t\right| ^{2H}\right) E\left(
X_{1}^{2}\right)  \label{2.10}
\end{equation}
The simplest process with these properties is a Gaussian process called
fractional Brownian motion. Fractional Brownian motion \cite{Mandelbrot1} 
\begin{equation}
\Bbb{E}\left[ B_{H}\left( t\right) \right] =0\qquad \Bbb{E}\left[
B_{H}\left( t\right) B_{H}\left( s\right) \right] =\frac{1}{2}\left\{ \left|
t\right| ^{2H}+\left| s\right| ^{2H}-\left| t-s\right| ^{2H}\right\}
\label{2.11}
\end{equation}
has, for $H>\frac{1}{2}$ , a long range dependence 
\begin{equation}
\sum_{n=1}^{\infty }\mathnormal{Cov}\left( B_{H}\left( 1\right) ,B_{H}\left(
n+1\right) -B_{H}\left( n\right) \right) =\infty  \label{2.12}
\end{equation}

(b) Therefore, mathematical simplicity suggests the identification of the $%
R_{\sigma }\left( t\right) $ process with fractional Brownian motion. 
\begin{equation}
R_{\sigma }\left( t\right) =kB_{H}\left( t\right)  \label{2.13}
\end{equation}
From the data one obtains the Hurst coefficient $H\simeq 0.8$ (for the NYSE
index). The same parametrization holds for the data of all individual
companies that were tested, with $H$ in the range $0.8-0.9$.

For comparison the plot in the down left panel of Fig.1 shows the scaling
test for the (NYSE) price process where, unlike the $R_{\sigma }\left(
t\right) $ process, clear deviations are seen on the first few days.

From (\ref{2.8}) and the identification (\ref{2.13}) one concludes that the
induced volatility may be modeled by 
\begin{equation}
\log \sigma _{t}=\beta +\frac{k}{\delta }\left( B_{H}\left( t\right)
-B_{H}\left( t-\delta \right) \right)  \label{2.14}
\end{equation}
$\delta $ being the observation time scale (one day, for daily data). It
means that the volatility is not driven by fractional Brownian motion but by
fractional noise. For the volatility (at resolution $\delta $) 
\begin{equation}
\sigma \left( t\right) =\theta e^{\frac{k}{\delta }\left\{ B_{H}\left(
t\right) -B_{H}\left( t-\delta \right) \right\} -\frac{1}{2}\left( \frac{k}{%
\delta }\right) ^{2}\delta ^{2H}}  \label{2.15}
\end{equation}
the term $-\frac{1}{2}\left( \frac{k}{\delta }\right) ^{2}\delta ^{2H}$
being included to insure that $E\left( \sigma \left( t\right) \right)
=\theta $.

Eqs. (\ref{2.1}) and (\ref{2.14}) define a stochastic volatility model.

\begin{equation}
\begin{array}{lll}
dS_{t} & = & \mu S_{t}dt+\sigma _{t}S_{t}dB\left( t\right) \\ 
\log \sigma _{t} & = & \beta +\frac{k}{\delta }\left\{ B_{H}\left( t\right)
-B_{H}\left( t-\delta \right) \right\}
\end{array}
\label{2.16}
\end{equation}
In this coupled stochastic system, in addition to a mean value, volatility
is driven by fractional noise. Notice that this empirically based model is
different from the usual stochastic volatility models which assume the
volatility to follow an arithmetic or geometric Brownian process. Also in
the Comte and Renault model\cite{Comte}, it is fractional Brownian motion
that drives the volatility, not its derivative (fractional noise). $\delta $
is the observation scale of the process. In the $\delta \rightarrow 0$ limit
the driving process would be the distribution-valued process $W_{H}$%
\begin{equation}
W_{H}=\lim_{\delta \rightarrow 0}\frac{1}{\delta }\left( B_{H}\left(
t\right) -B_{H}\left( t-\delta \right) \right)  \label{2.17}
\end{equation}

In (\ref{2.16}) the constant $k$ measures the strength of the volatility
randomness. Although phenomenologically grounded and mathematically well
specified, the stochastic system (\ref{2.16}) is still a limited model
because, in particular, the fact that the volatility is not correlated with
the price process excludes the modeling of leverage effects. It would be
simple to introduce, by hand, such a correlation in the second equation in (%
\ref{2.16}). However we do prefer not to do so at this time, because have
not yet found a natural way to do it, which is as clear-cut and imposed by
the data as the approach that led to (\ref{2.16}).

\section{The statistics of price returns}

Here one computes the probability distribution of price returns implied by
the stochastic volatility model (\ref{2.16}). From (\ref{2.14}) one
concludes that $\log \sigma _{t}$ is a Gaussian process with mean $\beta $
and covariance 
\begin{equation}
\psi \left( s,u\right) =\frac{k^{2}}{2\delta ^{2}}\left\{ \left| s-u+\delta
\right| ^{2H}+\left| u-s+\delta \right| ^{2H}-2\left| s-u\right|
^{2H}\right\}  \label{3.1}
\end{equation}
This Gaussian process has non-trivial correlation for $H\neq \frac{1}{2}$.
At each fixed time $\log \sigma _{t}$ is a Gaussian random variable with
mean $\beta $ and variance $k^{2}\delta ^{2H-2}$. Then, 
\begin{equation}
p_{\delta }\left( \sigma \right) =\frac{1}{\sigma }p_{\delta }\left( \log
\sigma \right) =\frac{1}{\sqrt{2\pi }\sigma k\delta ^{H-1}}\exp \left\{ -%
\frac{\left( \log \sigma -\beta \right) ^{2}}{2k^{2}\delta ^{2H-2}}\right\}
\label{3.2}
\end{equation}
therefore

\begin{equation}
P_{\delta }\left( \log \frac{S_{T}}{S_{t}}\right) =\int_{0}^{\infty }d\sigma
p_{\delta }\left( \sigma \right) p_{\sigma }\left( \log \frac{S_{T}}{S_{t}}%
\right)  \label{3.3}
\end{equation}
with 
\begin{equation}
p_{\sigma }\left( \log \frac{S_{T}}{S_{t}}\right) =\frac{1}{\sqrt{2\pi
\sigma ^{2}\left( T-t\right) }}\exp \left\{ -\frac{\left( \log \left( \frac{%
S_{T}}{S_{t}}\right) -\left( \mu -\frac{\sigma ^{2}}{2}\right) \left(
T-t\right) \right) ^{2}}{2\sigma ^{2}\left( T-t\right) }\right\}  \label{3.4}
\end{equation}
One sees that the effective probability distribution of the returns might
depend both on the time lag $\Delta =T-t$ and on the observation time scale $%
\delta $ used to construct the volatility process. That this latter
dependence might actually be very weak, seems to be implied by some
surprising experimental results.

\begin{figure}[tbh]
\begin{center}
\psfig{figure=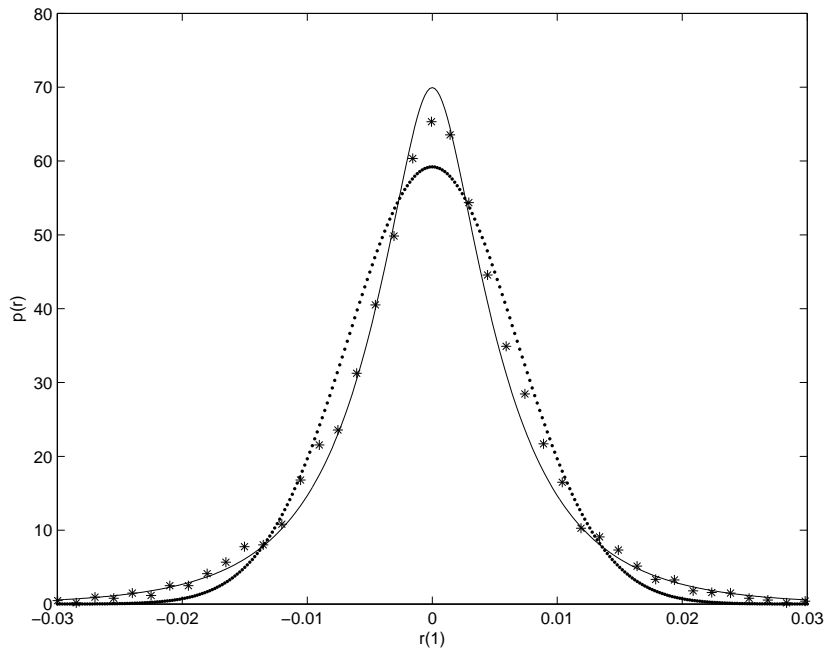,width=11truecm}
\end{center}
\caption{One-day NYSE returns compared with the model predictions and the
lognormal}
\end{figure}

Before obtaining a closed form expression for $P_{\sigma }\left( \log \frac{%
S_{T}}{S_{t}}\right) $ and its asymptotic behavior, we will present some
comparisons with market data. For the Fig.2 the same NYSE\ one-day data as
before is used to fix the parameters of the volatility process. Then, using $%
H=0.83,$ $k=0.59,$ $\beta =-5,$ $\delta =1$, the one-day return distribution
predicted by the model is compared with the data. The agreement is quite
reasonable. For comparison a log-normal with the same mean and variance is
also plotted in Fig.2. Then, in Fig. 3, using the same parameters, the same
comparison is made for the $\Delta =1$ and $\Delta =10$ data.

\begin{figure}[tbh]
\begin{center}
\psfig{figure=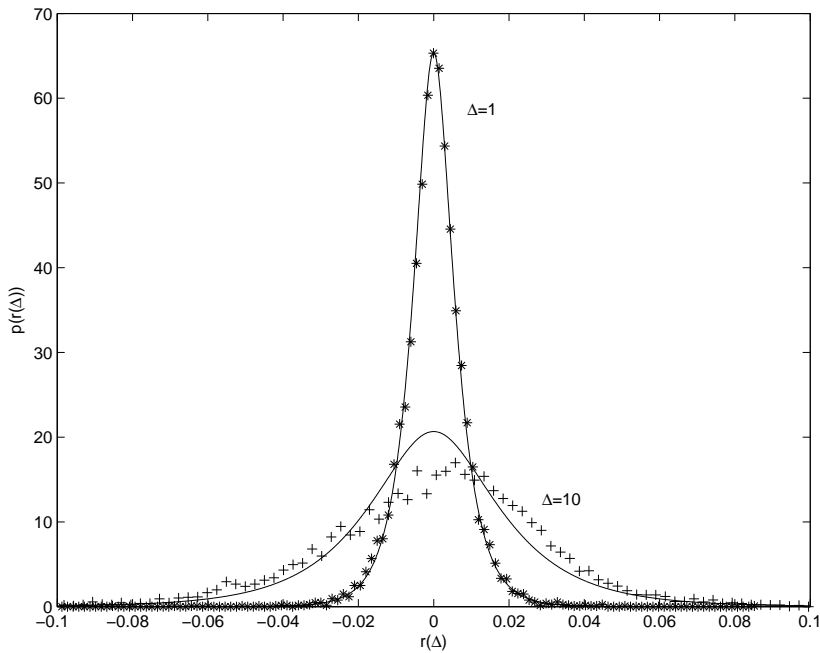,width=11truecm}
\end{center}
\caption{One and ten-days NYSE returns compared with the model predictions}
\end{figure}

\begin{figure}[tbh]
\begin{center}
\psfig{figure=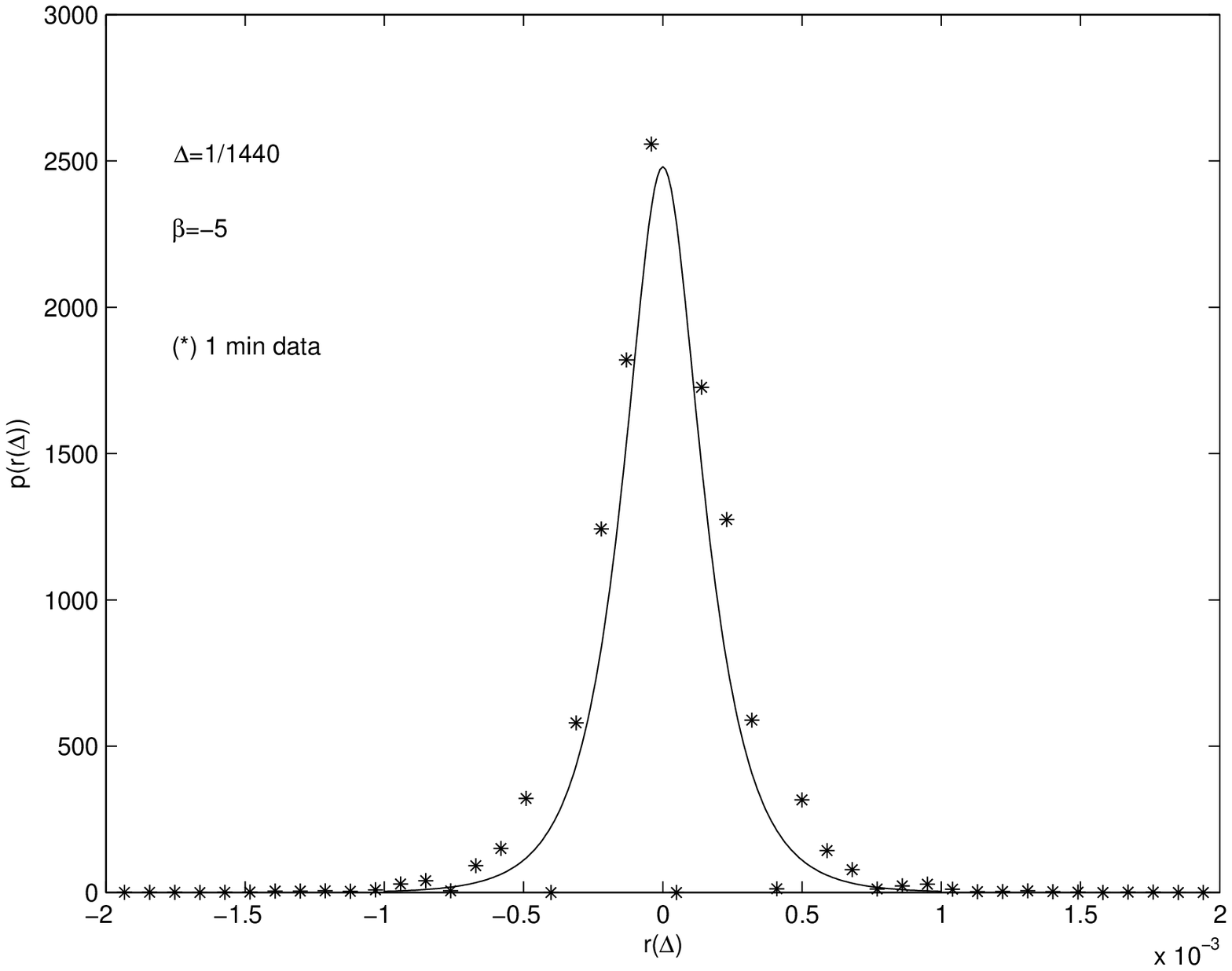,width=11truecm}
\end{center}
\caption{One-minute USD-Euro returns compared with the model predictions,
with parameters obtained from one-day NYSE data}
\end{figure}

Fig. 4 shows a somewhat surprising result. Using the same parameters and
just changing $\Delta $ from $1$ (one day) to $\Delta =\frac{1}{440}$ (one
minute), the prediction of the model is compared with one-minute data of
USDollar-Euro market for a couple of months in 2001. The result is
surprising, because one would not expect the volatility parametrization to
carry over to such a different time scale and also because one is dealing
with different markets. A systematic analysis of high-frequency data is now
being carried out to test the degree of time-scale dependence of the
volatility parametrization and its universality over different markets.

\begin{figure}[tbh]
\begin{center}
\psfig{figure=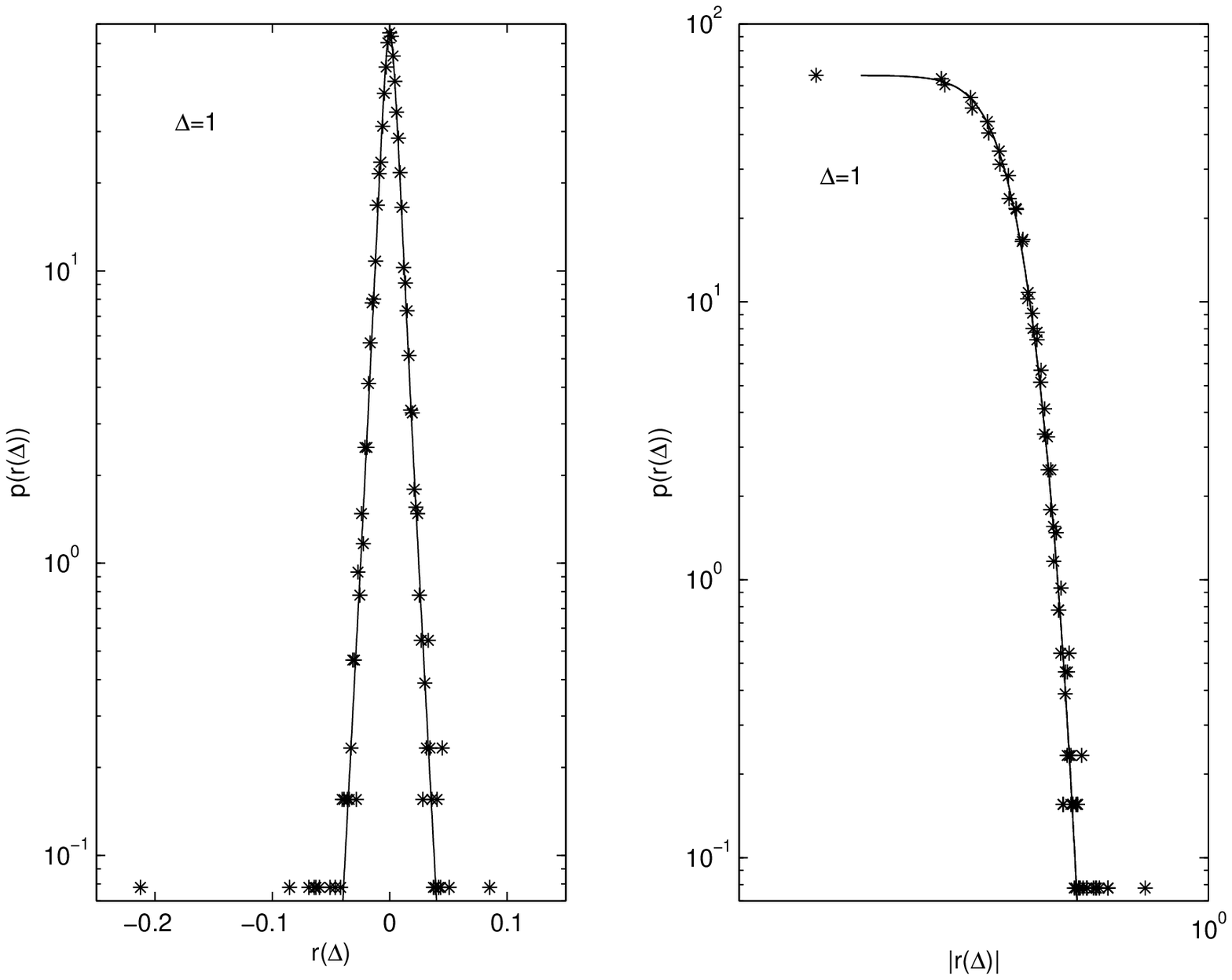,width=11truecm}
\end{center}
\caption{Semilogarithmic and loglog plots of NYSE data}
\end{figure}

\begin{figure}[tbh]
\begin{center}
\psfig{figure=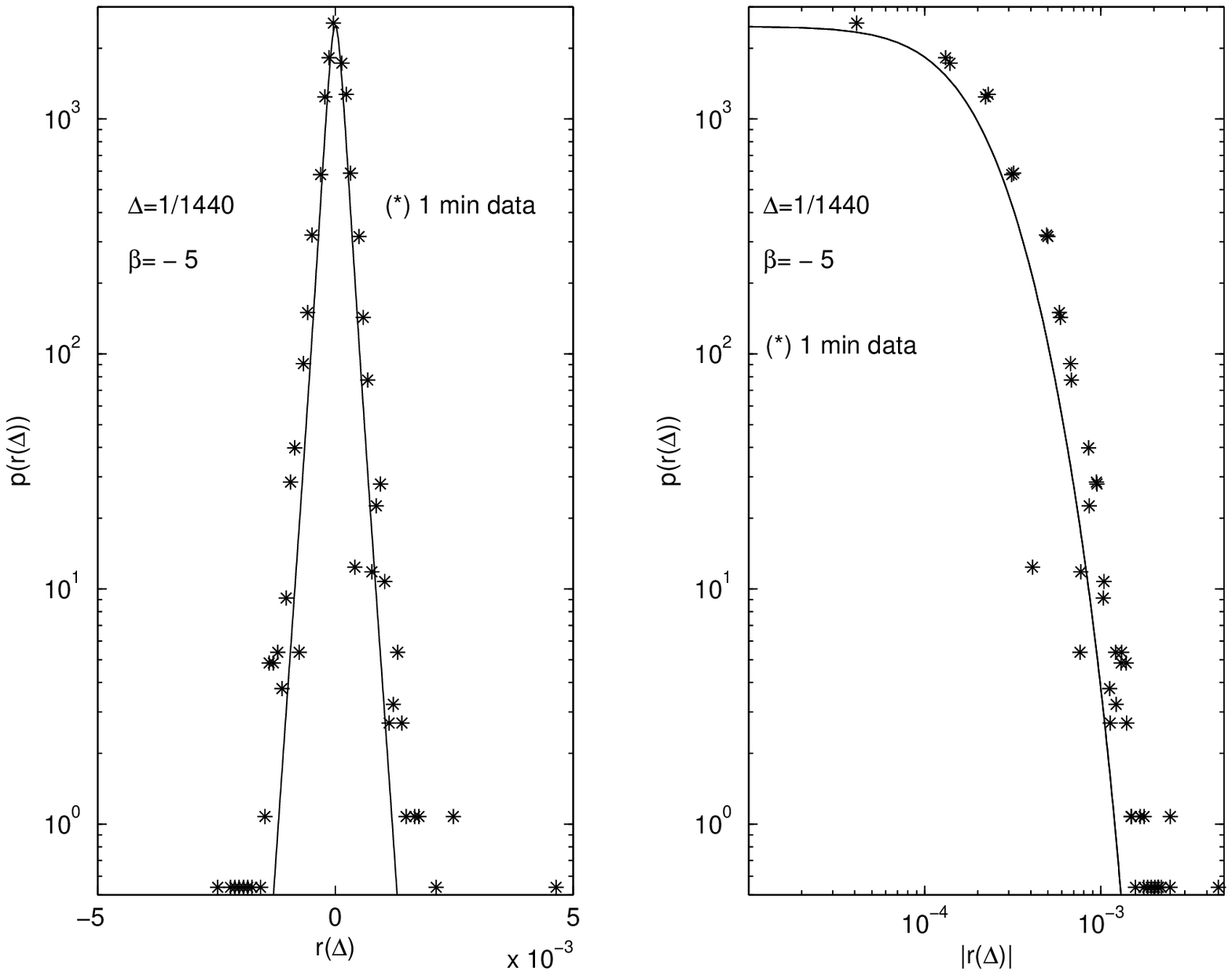,width=11truecm}
\end{center}
\caption{Semilogarithmic and loglog plots of USD-Euro data}
\end{figure}

In Fig.5 and Fig.6 we have displayed the same one-day and one-minute return
data discussed before as well as the predictions of the model both in
semilogarithmic and loglog plots.

Now we will establish a closed-form expression for the returns distribution
and its asymptotic behavior. Using (\ref{3.2}) and (\ref{3.4}) in (\ref{3.3}%
) and changing variables one obtains 
\begin{equation}
P_{\delta }\left( r\left( \Delta \right) \right) =\frac{1}{4\pi \theta
k\delta ^{H-1}\sqrt{\Delta }}\int_{0}^{\infty }dxx^{-\frac{1}{2}}e^{-\frac{1%
}{C}\left( \log x\right) ^{2}}e^{-\lambda x}  \label{3.10}
\end{equation}
with 
\begin{equation}
r\left( \Delta \right) =\log S_{T}-\log S_{t}\;,\;\theta =e^{\beta
}\;,\;\Delta =T-t\;,\;\lambda =\frac{\left( r\left( \Delta \right)
-r_{0}\right) ^{2}}{2\Delta \theta ^{2}}  \label{3.12}
\end{equation}
and 
\begin{equation}
r_{0}=\left( \mu -\frac{\sigma ^{2}}{2}\right) \left( T-t\right)
\;,\;C=8k^{2}\delta ^{2H-2}  \label{3.14}
\end{equation}
Expanding the exponential in (\ref{3.10}) 
\begin{eqnarray}
4\pi \theta k\delta ^{H-1}\sqrt{\Delta }P_{\delta }\left( r\left( \Delta
\right) \right) &=&\sum_{n=0}^{\infty }\frac{1}{n!}\left( -\frac{1}{C}%
\right) ^{n}\int_{0}^{\infty }dxx^{-\frac{1}{2}}e^{-\lambda x}\left( \log
x\right) ^{2n}  \nonumber \\
&=&\sum_{n=0}^{\infty }\frac{1}{n!}\left( -\frac{1}{C}\right) ^{n}\frac{%
\partial ^{2n}}{\partial z^{2n}}\left. \left( \lambda ^{-z}\Gamma \left(
z\right) \right) \right| _{z=\frac{1}{2}}  \label{3.14a}
\end{eqnarray}
Finally 
\begin{equation}
P_{\delta }\left( r\left( \Delta \right) \right) =\frac{1}{4\pi \theta
k\delta ^{H-1}\sqrt{\Delta }}\frac{1}{\sqrt{\lambda }}\left. \left( e^{-%
\frac{1}{C}\left( \log \lambda -\frac{d}{dz}\right) ^{2}}\Gamma \left(
z\right) \right) \right| _{z=\frac{1}{2}}  \label{3.15}
\end{equation}
with asymptotic behavior, for large returns 
\begin{equation}
P_{\delta }\left( r\left( \Delta \right) \right) \sim \frac{1}{\sqrt{\Delta
\lambda }}e^{-\frac{1}{C}\log ^{2}\lambda }  \label{3.16}
\end{equation}

On the other hand, as seen from Figs. 5 and 6, the exact result (\ref{3.10})
or (\ref{3.15}) resembles the double exponential distribution recognized by
Silva, Prange and Yakovenko\cite{Silva} as a new stylized fact in market
data. The double exponential distribution has been shown, by Dragulescu and
Yakovenko\cite{Dragulescu}, to follow from Heston's \cite{Heston} stochastic
volatility model. Notice however that our model is different from Heston's
model in that volatility is driven by a process with memory (fractional
noise). As a result, despite the qualitative similarity of behavior at
intermediate return ranges, the analytic form of the distribution and the
asymptotic behavior are different.

\section{Option pricing}

Assuming risk neutrality \cite{Cox1}, the value $V\left( S_{t},\sigma
_{t},t\right) $ of an option is the present value of the expected terminal
value discounted at the risk-free rate 
\begin{equation}
V\left( S_{t},\sigma _{t},t\right) =e^{-r\left( T-t\right) }\int V\left(
S_{T},\sigma _{T},T\right) p\left( S_{T}|S_{t},\sigma _{t}\right) dS_{T}
\label{4.1}
\end{equation}
$V\left( S_{T},\sigma _{T},T\right) =\max \left[ 0,S-K\right] $ and the
conditional probability for the terminal price depends on $S_{t}$ and $%
\sigma _{t}$. $K$ is the strike price, $T$ the maturity time and $S_{t}$ and 
$\sigma _{t}$ the price and volatility of the underlying security.

Whenever the drift of a financial time series can be replaced by the
risk-free rate we are in a risk-neutral situation. In stochastic volatility
models (with or without fractional noise) this is not an accurate
assumption. Nevertheless we will make use of (\ref{4.1}) to obtain an
approximate estimate of the deviations from Black-Scholes implied by the
stochastic differential model (\ref{2.16}). As in Hull and White \cite{Hull}%
, we make use of the relation between conditional probabilities of related
variables, namely 
\begin{equation}
p\left( S_{T}|S_{t},\sigma _{t}\right) =\int p\left( S_{T}|S_{t},\overline{%
\log \sigma }\right) p\left( \overline{\log \sigma }|\log \sigma _{t}\right)
d\left( \overline{\log \sigma }\right)  \label{4.2}
\end{equation}
$\overline{\log \sigma }$ being the random variable 
\begin{equation}
\overline{\log \sigma }=\frac{1}{T-t}\int_{t}^{T}\log \sigma _{s}ds
\label{4.3}
\end{equation}
that is, $\overline{\log \sigma }$ is the mean volatility from time $t$ to
the maturity time $T$ conditioned to an average value $\log \sigma _{t}$ at
time $t$. Then Eq.(\ref{4.1}) becomes 
\begin{equation}
V\left( S_{t},\sigma _{t},t\right) =\int C\left( S_{t},e^{\overline{\log
\sigma }},t\right) p\left( \overline{\log \sigma }|\log \sigma _{t}\right)
d\left( \overline{\log \sigma }\right)  \label{4.4}
\end{equation}
\begin{equation}
C\left( S_{t},e^{\overline{\log \sigma }},t\right) =\int e^{-r\left(
T-t\right) }V\left( S_{T},\sigma _{T},T\right) p\left( S_{T}|S_{t},\overline{%
\log \sigma }\right) dS_{T}  \label{4.5}
\end{equation}
$C\left( S_{t},e^{\overline{\log \sigma }},t\right) $ being the
Black-Scholes price for an option with average volatility $e^{\overline{\log
\sigma }}$, which is known to be \cite{Black} \cite{Merton} 
\begin{equation}
C\left( S_{t},\sigma ,t\right) =S_{t}\left( a+b\right) N\left( a,b\right)
-Ke^{-r\left( T-t\right) }\left( a-b\right) N\left( a,-b\right)  \label{4.6}
\end{equation}
with 
\begin{equation}
\begin{array}{lll}
a & = & \frac{1}{\sigma }\left( \frac{\log \frac{S_{t}}{K}}{\sqrt{T-t}}+r%
\sqrt{T-t}\right) \\ 
b & = & \frac{\sigma }{2}\sqrt{T-t}
\end{array}
\label{4.7}
\end{equation}
and 
\begin{equation}
N\left( a,b\right) =\frac{1}{\sqrt{2\pi }}\int_{-1}^{\infty }dye^{-\frac{%
y^{2}}{2}\left( a+b\right) ^{2}}  \label{4.8}
\end{equation}

In a stochastic volatility model with fractional noise, instead of $V\left(
S_{t},\sigma _{t},t\right) $, it would be more correct to write $V\left(
S_{t},\sigma _{\leq t},t\right) $ to emphasize the dependence on the past.
For simplicity we have used the first notation, with the provision that at
no point, in the calculation below, Markov properties of the processes
should be assumed, only their Gaussian nature.

To compute the conditional probability $p\left( \overline{\log \sigma }|\log
\sigma _{t}\right) $ it follows from (\ref{2.16}) that the process $%
\overline{\log \sigma }$ conditioned to $\log \sigma _{t}$ at $t$ is 
\begin{equation}
\overline{\log \sigma }=\log \sigma _{t}+\frac{1}{T-t}\int_{t}^{T}\frac{k}{%
\delta }ds\int_{t}^{s}\left( dB_{H}\left( \tau \right) -dB_{H}\left( \tau
-\delta \right) \right)  \label{4.9}
\end{equation}
Notice that, because we want to compute the conditional probability of $%
\overline{\log \sigma }$ given $\log \sigma _{t}$ at time $t$, $\sigma _{t}$
in Eq.(\ref{4.9}) is not a process but simply the value of the argument in
the $V\left( S_{t},\sigma _{t},t\right) $ function.

As a $t-$dependent process the double integral in (\ref{4.9}) is a centered
Gaussian process. Therefore, given $\log \sigma _{t}$ at time $t$, $%
\overline{\log \sigma }$ is a Gaussian variable with conditional mean and
variance 
\begin{equation}
E\left\{ \overline{\log \sigma }|\log \sigma _{t}\right\} =\log \sigma _{t}
\label{4.10}
\end{equation}
\begin{eqnarray}
\alpha ^{2} &=&E\left\{ \left( \overline{\log \sigma }-\log \sigma
_{t}\right) ^{2}|\log \sigma _{t}\right\}  \nonumber \\
&=&\frac{k^{2}}{\delta ^{2}\left( T-t\right) ^{2}}E\left\{
\int_{t}^{T}ds\int_{t}^{s}\left[ dB_{H}\left( \tau \right) -dB_{H}\left(
\tau -\delta \right) \right] \int_{t}^{T}ds^{^{\prime
}}\int_{t}^{s^{^{\prime }}}\left[ dB_{H}\left( \tau ^{^{\prime }}\right)
-dB_{H}\left( \tau ^{^{\prime }}-\delta \right) \right] \right\}
\label{4.10a}
\end{eqnarray}
Expanding $\int_{t}^{s}\left[ dB_{H}\left( \tau \right) -dB_{H}\left( \tau
-\delta \right) \right] =B_{H}\left( s\right) -B_{H}\left( t\right)
-B_{H}\left( s-\delta \right) +B_{H}\left( t-\delta \right) $ and using (\ref
{2.11}) one obtains 
\begin{equation}
\alpha ^{2}=\frac{k^{2}}{\delta ^{2}\left( T-t\right) }\left\{ \frac{1}{%
2\left( T-t\right) }I_{1}+I_{2}\right\} +k^{2}\delta ^{2H-2}  \label{4.11}
\end{equation}
with 
\begin{equation}
I_{1}=\frac{2}{\left( 2H+1\right) \left( 2H+2\right) }\left\{ \left(
T-t+\delta \right) ^{2H+2}+\left( T-t-\delta \right) ^{2H+2}-2\left(
T-t\right) ^{2H+2}-2\delta ^{2H+2}\right\}  \label{4.12}
\end{equation}
\[
I_{2}=\frac{1}{2H+1}\left\{ 2\left( T-t\right) ^{2H+1}-\left( T-t+\delta
\right) ^{2H+1}-\left( T-t-\delta \right) ^{2H+1}\right\} 
\]
As is seen by expanding $I_{1}$ and $I_{2}$, when $t\rightarrow T$ one has
the consistency condition $\alpha ^{2}\rightarrow 0$. However, in general,
for option pricing purposes, $\delta \ll \left( T-t\right) $ and one may
approximate 
\begin{equation}
\alpha ^{2}\simeq \frac{k^{2}}{\delta ^{2-2H}}\left( 1-\left( 2H-1\right)
\left( \frac{\delta }{T-t}\right) ^{2-2H}\right)  \label{4.13}
\end{equation}
Finally 
\begin{equation}
p\left( \overline{\log \sigma }|\log \sigma _{t}\right) =\frac{1}{\sqrt{2\pi 
}\alpha }\exp \left\{ \frac{-\left( \overline{\log \sigma }-\log \sigma
_{t}\right) ^{2}}{2\alpha ^{2}}\right\}  \label{4.14}
\end{equation}
and from (\ref{4.4}) 
\begin{equation}
V\left( S_{t},\sigma _{t},t\right) =\int_{-\infty }^{\infty }d\xi C\left(
S_{t},e^{\xi },t\right) p\left( \xi |\log \sigma _{t}\right)  \label{4.15}
\end{equation}
one obtains 
\begin{equation}
V\left( S_{t},\sigma _{t},t\right) =S_{t}\left[ aM\left( \alpha ,a,b\right)
+bM\left( \alpha ,b,a\right) \right] -Ke^{-r\left( T-t\right) }\left[
aM\left( \alpha ,a,-b\right) -bM\left( \alpha ,-b,a\right) \right]
\label{4.16}
\end{equation}
\begin{eqnarray}
M\left( \alpha ,a,b\right) &=&\frac{1}{2\pi \alpha }\int_{-1}^{\infty
}dy\int_{0}^{\infty }dxe^{-\frac{\log ^{2}x}{2\alpha ^{2}}}e^{-\frac{y^{2}}{2%
}\left( ax+\frac{b}{x}\right) ^{2}}  \label{4.17} \\
&=&\frac{1}{4\alpha }\sqrt{\frac{2}{\pi }}\int_{0}^{\infty }dx\frac{e^{-%
\frac{\log ^{2}x}{2\alpha ^{2}}}}{ax+\frac{b}{x}}\mathnormal{\mathnormal{erf}%
\mathnormal{c}}\left( -\frac{ax}{\sqrt{2}}-\frac{b}{\sqrt{2}x}\right) 
\nonumber
\end{eqnarray}
as a new option price formula (\textnormal{erf}c is the complementary error
function and $a$ and $b$ are defined in Eq.(\ref{4.7})), with $\sigma $
replaced by $\sigma _{t}$.

Eqs. (\ref{4.15}) and (\ref{4.16}) are mathematically equivalent. For
computational convenience (of the reader that might want to use our formula)
we point out that, instead of writting performing codes for the M-functions
in Eq.(\ref{4.16}), he might simply use a Black-Scholes code and perform the
integration in Eq.(\ref{4.15}).

\begin{figure}[htb]
\begin{center}
\psfig{figure=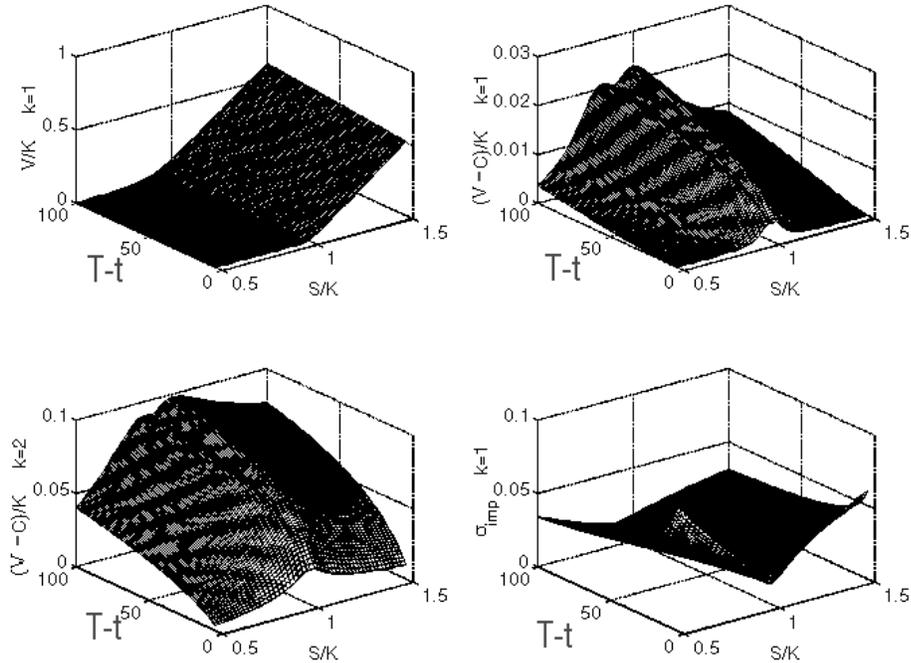,width=11truecm}
\end{center}
\caption{Option price and equivalent implied volatility in the
"risk-neutral" approach to the stochastic volatility model}
\end{figure}

In Fig.7 we plot the option value surface for $V\left( S_{t},\sigma
_{t},t\right) $ in the range $T-t\in [5,100]$ and $S/K\in [0.5,1.5]$ as well
as the difference $\left( V\left( S_{t},\sigma _{t},t\right) -C\left(
S_{t},\sigma _{t},t\right) \right) /K$ for $k=1$ and $k=2$. The other
parameters are fixed at $\sigma =0.01,r=0.001,\delta =1,H=0.8$.

To compare the predictions of our formula with the classical Black-Scholes
(BS) result, we have computed the implied volatility required in the BS
model to reproduce our results. This is plotted in the lower right panel of
Fig.7 which shows the implied volatility surface corresponding to $V\left(
S_{t},\sigma _{t},t\right) $ for $k=1$. One sees that, when compared to BS,
it predicts a smile effect with the smile increasing as maturity approaches.

\section{Conclusions}

(a) In this paper, rather than starting by postulating some model for the
market process and then exploring its better or worse vindication by the
data, the approach has been to be inspired, at each step of its
construction, both by mathematical simplicity and consistency with the data.
It is mathematically more complex and requires (for example for a derivation
of option pricing without assuming risk-neutrality) more sophisticated tools
of Malliavin calculus than most stochastic volatility models. Nevertheless,
from its very construction and consistency with the data, it appears as a
kind of minimal model.

(b) The asymptotic behavior of price returns, in special its asymptotic
behavior has been much discussed (see for example \cite{Toyli} and
references therein). In particular it has been proposed that the large
return tail decays as a power law, although a stretched exponential might
provide a better fit\cite{Malevergne}.

The semilogarithmic plots suggest in fact that a better overall fit might be
obtained by a stretched exponential or indeed by Eqs. (\ref{3.15}) and (\ref
{3.16}).

(c) From the data and model comparison plotted in the figures it looks that
the stochastic volatility model (as well as a scaling hypothesis) cannot fit
the very large deviations. There is a good fit for the bulk of the data but
there are also a few events very far from the fit. It suggests that a model
with two probability spaces is still not enough to capture the whole
process. Maybe one should write $S_{t}\left( \omega ,\omega ^{^{\prime
}},\omega ^{^{\prime \prime }}\right) $ with the last entry, $\omega
^{^{\prime \prime }}$, representing exogenous market shocks.

\end{document}